\documentclass[10pt]{article}

\usepackage{amsmath}
\usepackage{amssymb}
\usepackage{graphicx}

\usepackage{cite}

\usepackage{color} 
\usepackage{enumitem}
\usepackage{tikz}
\usepackage{srcltx}

\usepackage[labelfont=bf,labelsep=period,justification=raggedright]{caption}

\makeatletter
\renewcommand{\@biblabel}[1]{\quad#1.}
\makeatother

\date{}

\newcommand{\abs}[1]{\lvert#1\rvert}

\newcommand{\sgn}{\operatorname{sgn}}
\usepackage{esint,amssymb}
\newcommand{\Z}{\mathbb{Z}}
\newcommand{\N}{\mathbb{N}}

\begin{document}

% Title must be 150 characters or less
\begin{flushleft}
{\Large
\textbf{Traveling pulses for a two-species chemotaxis model}
}
% Insert Author names, affiliations and corresponding author email.
\\

\begin{center}
Casimir Emako$^{1}$, 
Charl\`ene Gayrard$^{2}$, 
Axel Buguin$^{3,\ast}$,
Lu\'\i s Neves de Almeida $^{4,\ast}$,
Nicolas Vauchelet $^{1,\ast}$
\end{center}

\textit{1} Sorbonne Universit\'es, UPMC Univ Paris 06, CNRS, INRIA,  UMR 7598, Laboratoire Jacques-Louis Lions, Equipe MAMBA, 4, place Jussieu 75005, Paris, France.
\\
\textit{2} Institut Jacques Monod - UMR 7592, CNRS – Universit\'e Denis Diderot, Paris, France
\\
\textit{3} Laboratoire Physico-chimie Curie, Institut Curie, PSL Research University, CNRS UMR 168, 75005 Paris, France
\\
\textit{4} Sorbonne Universit\'es, CNRS, UPMC Univ Paris 06, INRIA,  UMR 7598, Laboratoire Jacques-Louis Lions, Equipe MAMBA, 4, place Jussieu 75005, Paris, France.
\smallskip \\

$\ast$ Senior authors in alphabetical order. Corresponding author email address: luis.almeida@upmc.fr
\end{flushleft}

% Please keep the abstract between 250 and 300 words
\section*{Abstract}
Mathematical models have been widely used to describe the collective movement of bacteria 
by chemotaxis. In particular, bacterial concentration waves traveling in a narrow channel have been experimentally 
observed and can be precisely described thanks to a mathematical model at the macroscopic scale.
Such model was derived in \cite{Saragosti} using a kinetic model based on an accurate description of the mesoscopic run-and-tumble process.
We extend this approach to study the behavior of the interaction between two populations of E. Coli.
Separately, each population travels with its own speed in the channel.
When put together, a synchronization of the speed of the traveling pulses can be observed.
We show that this synchronization depends on the fraction of the fast population.
Our approach is based on mathematical analysis of a macroscopic model of partial differential equations.
Numerical simulations in comparison with experimental observations show qualitative agreement.

% Please keep the Author Summary between 150 and 200 words
% Use first person. PLoS ONE authors please skip this step. 
% Author Summary not valid for PLoS ONE submissions.   

\section*{Author Summary}

The use of mathematical tools to describe self-organization of bacterial communities has raised 
a lot of interest since it permits a precise description of experimentally observed phenomena. 
In the last 40 years a hierarchy of mathematical models for the dynamics of a single bacterial population has been proposed. 
These models have progressively taken into account more precise aspects of individual bacterial behavior (like the run and tumble behavior). 
Nowadays, a natural and challenging issue is to use such models 
to describe the interaction between different populations of bacteria. 
In this work, we consider a macroscopic system of equations derived from the mesoscopic scales 
to describe the interaction between two populations of bacteria.
The prediction obtained thanks to this model is compared to experimental observations
concerning the behavior of traveling pulses of bacteria in a channel.

\newpage

\section*{Introduction}
\setcounter{section}{1}
\setcounter{subsection}{0}
The ability of microorganisms to sense their environment helps them to colonize regions by using chemical cues to move towards favorable areas (e.g. with higher concentration in nutrients like in the present study).
This biological process called chemotaxis has been extensively studied. Since pioneering works of Adler \cite{Adler} (see also \cite{BudBerg, Keller}), 
we know that many bacteria like E.Coli may gather, feel the nutrient (oxygen, glucose) and move towards it. 
Many models including the famous Keller-Segel system (see \cite{Keller-Segel,Hillen,Tindall,Benoit}) were proposed to describe mathematically this behavior. 
In \cite{VCJS}, Saragosti et al. described the propagation of bacterial concentration waves in micro-channels using a macroscopic model. This model \eqref{one_species_model}
gives the dynamics of the density of cells $\rho(x,t)$ and the concentrations in chemoattractant $S(x,t)$ and nutrients~$N(x,t)$
\begin{equation}\label{one_species_model}
 \left\{
\begin{aligned}
\partial_t \rho&=D \Delta  \rho-\nabla \cdot \left(\rho \left(u[S]+u[N]\right)\right), \\
\partial_t S&=D_S\Delta S-\alpha S+\rho,\\
\partial_t N&=D_N\Delta N -\gamma \rho N,
\end{aligned}
\right.
\end{equation}
where $u[S],u[N]$ are given by
\begin{equation*}
 u[S]=\chi^S \sgn(\partial_x S),\quad u[N]=\chi^N \sgn(\partial_x N),
\end{equation*}
with $D,D_S,D_N,\chi^S,\chi^N,\alpha,\gamma$ positive constants. Here $\sgn$ is the sign function:
\begin{equation*}
\sgn(x)=
\left\{
 \begin{aligned}
  1,\quad &\text{if }x>0,\\
  -1,\quad &\text{if }x<0.
 \end{aligned}
\right.
\end{equation*}
The velocity fields $u[S]$ and $u[N]$ model the biased motion of bacteria
due to the attraction of chemoattractant and nutrient. The function $\sgn$ allows to model 
the fast response of bacteria to variations in their environment; this choice is motivated by the 
comparisons with experimental data presented in \cite{Saragosti,VCJS} which shows a good agreement with numerical simulations for this macroscopic model.

\begin{figure}[ht!]
\setlength\abovecaptionskip{-2ex}
 \centering
 \includegraphics[scale=0.5]{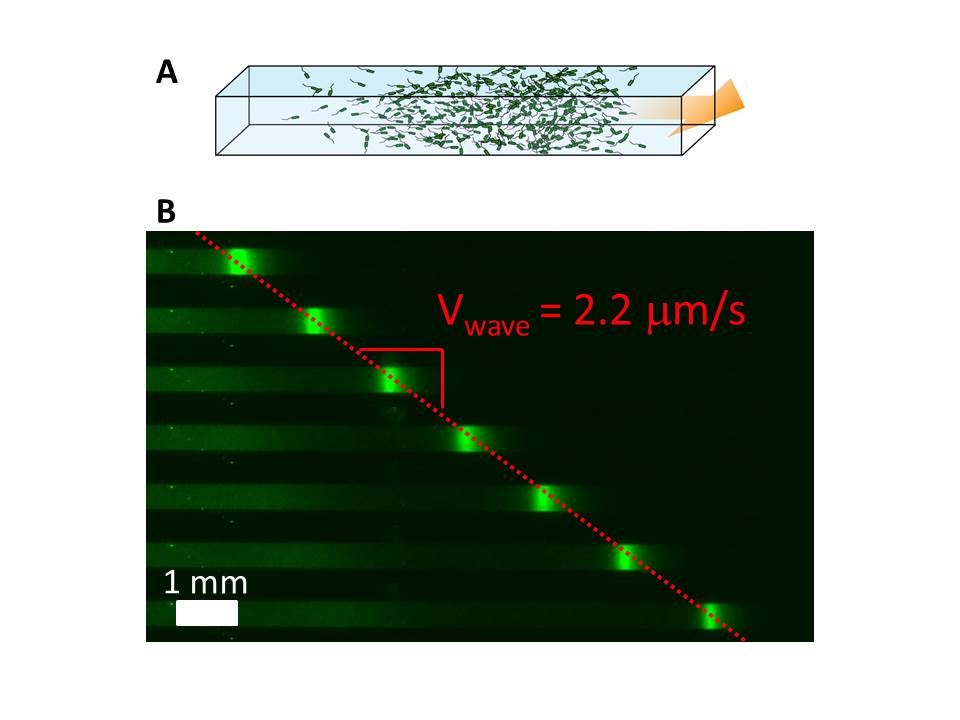}
 \caption{Collective migration of Escherichia coli in a PDMS micro-channel:\\
(A) Schematic view of a portion of the micro-channel. The channel is filled with a homogeneous suspension of bacteria and then centrifuged to accumulate bacteria on the left end. Few minutes after the centrifugation has been stopped, a concentration wave of bacteria propagates at constant velocity from left to right (see \cite{VCJS}).\\
(B) A concentration “wave” of E. coli chemotactic bacteria observed by fluorescence microscopy (white scale bar = 1 mm) propagating inside the micro-channel (top view). 
Successive snapshots of the same channel (600 s between successive images, the fluorescence intensity reflects the local concentration). The population migrates at a constant velocity (Vwave=2.2 $\mu$m/s for this particular experiment).\medskip\\ }
 \label{one_species}
\end{figure}

When initially put on the left of a channel filled with nutrients, the bacteria consume nutrients located at this side. This creates a gradient of nutrients (oriented) towards the right which 
induces the motion of bacteria represented by the drift term $u[N]$ in the equation for $\rho$ (see Fig \ref{one_species}).
While traveling in the channel, bacteria stay together thanks to the chemoattractant S they produce. 
From a mathematical point of view, particular solutions in translation with
constant velocity are called traveling waves. Such solutions have been studied for a long time
for reaction-diffusion equations, since the seminal work of R. A. Fisher \cite{Fisher} (see also \cite{KPP}). 
However, in such models, the motion of the front is induced
by the reaction term whereas in the case of model \eqref{one_species_model} there is no reaction term: the displacement by chemotaxis is modeled by the drift term in the equation for the bacteria density.
In \cite{VCJS}, existence of traveling waves for model \eqref{one_species_model} is proved, 
analytical forms for $\rho$ and $S$ are provided and the speed of the wave $\sigma$ 
satisfies the following relationship
\begin{equation}\label{sigma_one_species}
 \chi^N-\sigma=\chi^S \frac{\sigma}{\sqrt{4D_S\alpha +\sigma^2}}.
\end{equation}
This model predicts the double asymmetric exponential shape of $\rho$ and the speed $\sigma$ as observed experimentally.
However, when we have several different strains of bacteria simultaneously, new behaviors may emerge and this model might no longer be valid. 
For instance, if we start now from a mix of two bacterial strains having different swimming speeds what happens? Do they swim together? How does one population affects the other? This is the type of questions we want to address in the present paper.

\begin{figure}[ht!]
\centering
 \includegraphics[scale=0.5]{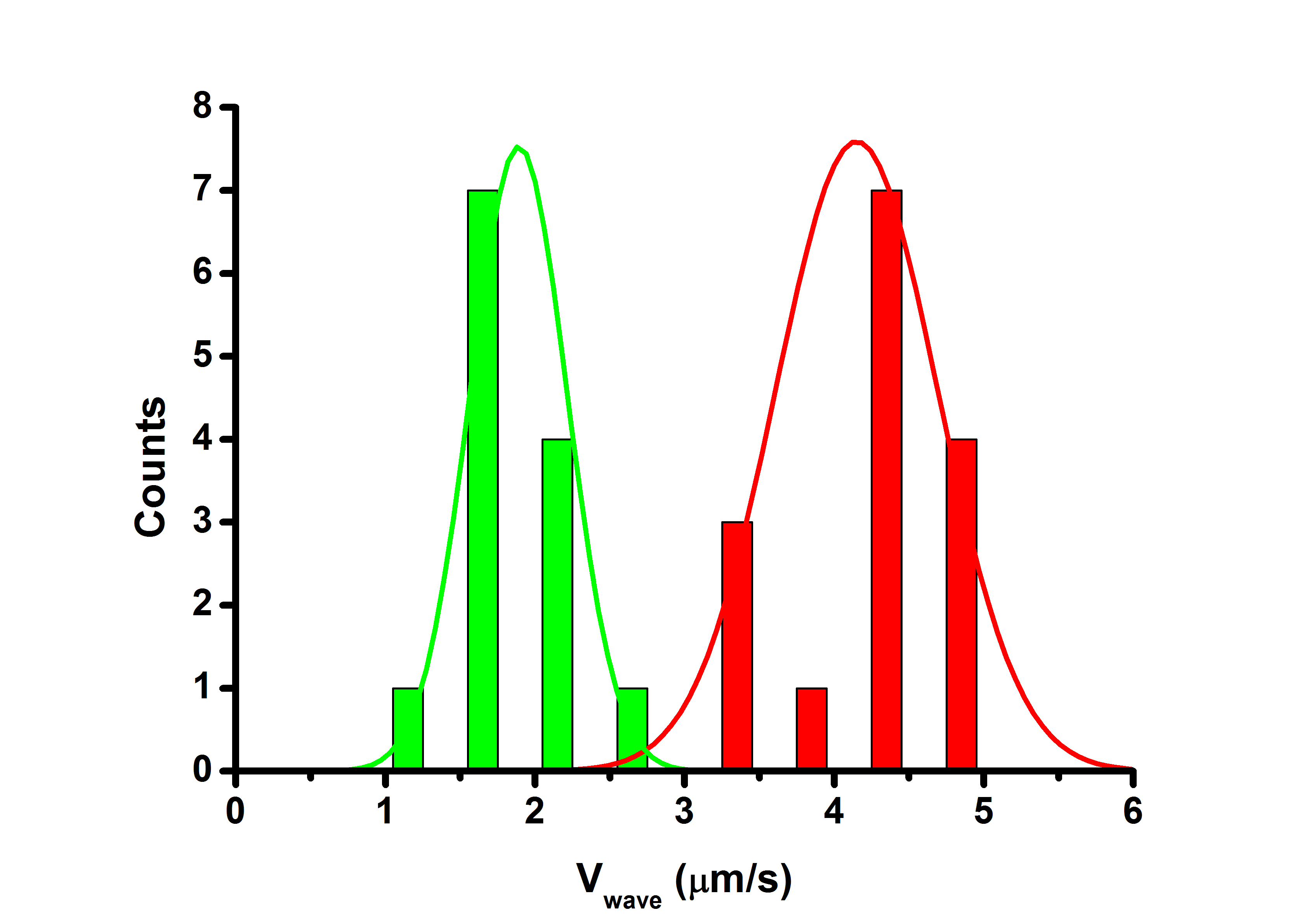}
 \caption{Wave velocity distribution for two different bacteria:
The same strain (RP437) transformed with two plasmids (PZE1R-GFP in green and PZE1R-mCherry in red) exhibits two different velocities for the wave propagation. 
The mean wave velocities obtained from about 15 experiments for each color, are respectively Vgreen=1.9 $\mu$m/s and Vred=4.1 $\mu$m/s.\medskip\\}
 \label{two_speeds}
\end{figure}

In this work, we study the case of a population composed of two different subpopulations of E.Coli that, when they are alone, form bands traveling at different speeds (subpopulation 1, green,  being the one traveling at a slower speed $\sigma_1$
and subpopulation 2, red,  being the one traveling at a higher speed $\sigma_2$, see Fig \ref{two_speeds}).\\

For the two strains used in this study, we observed that even if separately they travel at different speeds, when they are in presence of each other
they may form a single band and travel together. More precisely, our experiments show that when the ratio between the number of individuals of type 2 (fast subpopulation) and the number of individuals of type 1 (slow subpopulation) is sufficiently small, there is a single band
(i.e the two subpopulations travel together with an intermediate speed $\sigma$ such that $\sigma_1<\sigma<\sigma_2$).
We provide a mathematical model to describe this behavior and handle the dependency of the speed of the wave on relative sizes of subpopulations.

On the other hand, when this ratio is big, our experiments and our numerical simulations (Fig \ref{two_regimes}) show that the two subpopulations travel at different speeds.

\section*{Results}
\setcounter{section}{2}
\setcounter{subsection}{0}
\subsection{Description of the experiments}
When confined at one end of a micro-channel, large enough populations of swimming bacteria E. Coli propagate as concentration waves. 
To perform such experiments, we simply fill a micro-channel (in our experimental setting they have height = 100 $\mu$m, width = 500 $\mu$m and total length = 1.8 cm and are micro-fabricated using soft lithography \cite{Xia}) with a homogeneous solution of bacteria grown up to the mid-log phase ($5\times 10^8$ bacteria/mL). 
The channel is then closed at both ends using epoxy glue and gently centrifuged to accumulate motile bacteria at one end of the channel. When the centrifugation is stopped a concentration wave propagates along the channel at a velocity of a few $\mu$m/s. 
We use fluorescently labelled bacteria and thus it is possible to characterize the concentration profile of the traveling pulse using fluorescence video microscopy. These experiments are reported in previous publications \cite{Saragosti,VCJS}. 

In this work, in order to consider the case of multiple subpopulations, we used two types of bacteria: one carrying a plasmid expressing GFP (green) and the other carrying a plasmid expressing mCherry (red). The concentration waves obtained with the red ones are two times faster than the concentration waves of the green ones (see Fig \ref{two_speeds}). 
In the present paper we study, both experimentally and with our mathematical model, the behavior of the concentration waves obtained for different ratios $\phi_{red}=\frac{M_{red}}{M_{red}+M_{green}}$ (where $M_{red}$ and $M_{green}$ are the sizes of the two sub-population of bacteria) keeping constant the total number $(M_{red}+M_{green})$ of bacteria.

\begin{figure}[ht!]
\setlength\abovecaptionskip{-8ex}
\centering
 \includegraphics[scale=0.6]{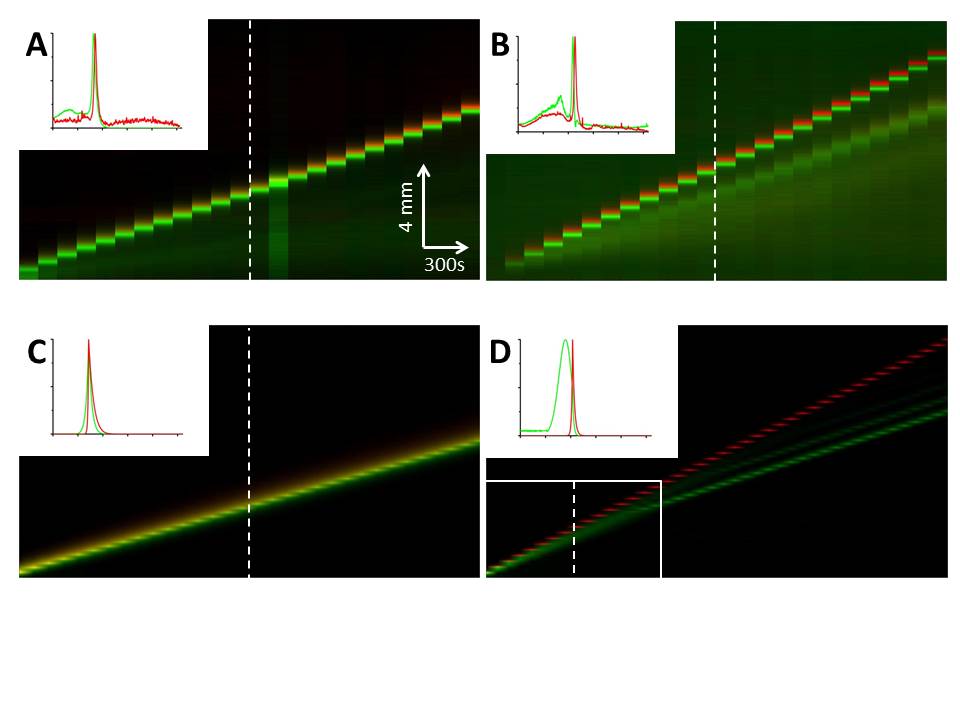}
 \caption{Kymographs showing the wave behavior for different bacterial compositions:
$\phi_{red}$ corresponds to the ratio of red bacteria over the total number of bacteria at the beginning of an experiment.\\
(A) Experimental result obtained with $\phi_{red}=10\%$.\\
(B) Experimental result obtained with $\phi_{red}=90\%$.\\
(C) Simulation based on our model with $\phi_{red}=10\%$.\\
(D) Simulation with $\phi_{red}=90\%$. 
The white rectangle in the lower left corner corresponds to the size of images A, B and C. In this case, the simulation is extended to show what happens at longer timescales (beyond the transitory regime so that the peaks are well separated). The small shift between the green and the red fronts in A and B (less visible in the insets) is due to the fact that while acquiring the images, we switch between the green and the red fluorescence channels every minute.\\
Insets: typical concentration profiles (corresponding to bacteria concentration profiles along the white dashed line for each kymograph). We normalized the green and the red peaks (so that the maximum value is one). The length of the the horizontal axis on the profile plots is 1.3 cm (this corresponds to the length of the dashed lines on the kymographs).
\\
}
\medskip
\label{two_regimes}
\end{figure}

\subsection{The one species model}

For the reader's convenience, we recall briefly in this subsection the main results obtained 
in \cite{VCJS}. That work establishes the existence of traveling waves for the one species model
\eqref{one_species_model}. 
Looking for traveling wave solutions to system \eqref{one_species_model} boils down 
to looking for particular solutions of the form $\rho(t,x)=\tilde{\rho}(x-\sigma t)$, 
$S(t,x)=\tilde{S}(x-\sigma t)$ and $N(t,x)=\tilde{N}(x-\sigma t)$. Moreover, since we are 
looking for a pulse, we have $\lim_{|z|\to \infty} \tilde{\rho}(z)=0$.
Injecting these expressions into the first equation of system \eqref{one_species_model}  
we get, after one integration,
$$
-\sigma \tilde{\rho} = D \tilde{\rho}' - \tilde{\rho} 
(\chi^S \sgn(\tilde{S}') +\chi^N \sgn(\tilde{N}')).
$$
To solve this equation, we make the ansatz that the wave moves from the left to the right, i.e. the gradient of the nutrient is positive ($\tilde{N}'>0$), and suppose that $\rho$ and $S$ 
are maximal at the same point (which, by translational invariance, is assumed to be $0$). 
Then, on $(-\infty,0)$, we have $\sgn(\tilde{S}')=\sgn(\tilde{N}')=1$, 
on $(0,+\infty )$, we have $\sgn(\tilde{S}')=-1$ and $\sgn(\tilde{N}')=1$.
Solving the differential equation, we have
\begin{equation}\label{exprho}
\tilde\rho(z)=e^{(\chi^N +\chi^S-\sigma)z/D}, \quad \mbox{ for } z<0; \qquad
\tilde\rho(z)=e^{(\chi^N -\chi^S-\sigma)z/D}, \quad \mbox{ for } z>0.
\end{equation}
In order to satisfy the vanishing condition at infinity, the velocity $\sigma$
should satisfy the condition $\chi^N-\chi^S < \sigma < \chi^N+\chi^S$.
Finally, we compute the velocity $\sigma$.
To do so, since $\tilde{S}$ is maximal at $0$, we have $\tilde{S}'(0)=0$.
From the expression for $\tilde{\rho}$ in \eqref{exprho}, we may solve the equation 
for $\tilde{S}$. Then, the condition $\tilde{S}'(0)=0$ gives a nonlinear problem 
for the velocity $\sigma$. After tedious but straightforward computations (see
the supplementary materials of \cite{VCJS}), we obtain and have to solve equation 
\eqref{sigma_one_species}; since its left hand side 
is decreasing with respect to $\sigma$ whereas its right hand side is nondecreasing, 
there exists a unique traveling speed $\sigma$ solving equation \eqref{sigma_one_species}.

Finally, we notice that we have a simple explicit expression of the density profile in
\eqref{exprho}. In \cite{VCJS} (in particular see Fig 2) this profile was compared to the one observed experimentally,  
showing a good agreement between experimental and analytical results.

\subsection{Description of the model}
We study the migration of a bacterial population composed of two subpopulations which react to two common chemical substances: the chemoattractant $S$ and the nutrient $N$.
These two chemical substances play different roles since bacteria produce the same chemoattractant which gathers the population and at the same time they consume the common nutrient which triggers the motion.
Each species is represented by its density at position $x \in \mathbb{R}^d$ and time $t>0$, $\rho_i(x,t)$ for $i=1,2$. The chemoattractant and the nutrient are described respectively by their concentration $S(x,t)$ and $N(x,t)$. 
Dynamics of $\rho_1,\rho_2,S,N$ are given by coupled advection-diffusion-reaction equations:
\begin{equation}\label{two_species_model}
\left\{
\begin{aligned}
\partial_t \rho_{1}&=D_1 \Delta  \rho_{1}-\nabla \cdot \left(\rho_{1}\left(u_1[S]+u_1[N]\right)\right), \\
\partial_t \rho_{2}&=D_2 \Delta  \rho_{2}-\nabla \cdot \left(\rho_{2}\left(u_2[S]+u_2[N]\right)\right),\\
\partial_t S&=D_S\Delta S-\alpha S+\rho_1+\rho_2,\\
\partial_t N&=D_N\Delta N -\gamma_1 \rho_1 N-\gamma_2 \rho_2 N,
\end{aligned}
\right.
\end{equation}
where $D_1,D_2,D_S,D_N,\alpha,\gamma_1,\gamma_2$ are positive constants.\\
The model is an extension to the two subpopulation case of the single population  model \eqref{one_species_model}. 
The starting point of the modeling is the Othmer-Dunbar-Alt model which illustrates the run-and-tumble process characterizing the motion of individual bacteria.
Then, the macroscopic equation is recovered by drift-diffusion limits.
By this means, we derive our model \eqref{two_species_model} from kinetic equations describing the phenomenon at the microscopic scale (see \cite{AEV}) and obtain expressions for $u_i[S],u_i[N]$  (this derivation is detailed in subsection \ref{sec:kin}).
Since the phenomenon we are considering is uni-directional, we restrict our study to one dimension in space ($d=1$) and for computational purposes,
we consider the following particular forms of $u_i[S]$ and $u_i[N]$
\begin{equation}\label{signals}
 u_i[S]=\chi_i^S\sgn (\partial_x S),\quad u_i[N]=\chi_i^N\sgn (\partial_x N),\qquad i=1,2,
\end{equation}
with $\chi_i^S,\chi_i^N$, the chemotactic sensitivities of subpopulation $i$ to the chemoattractant $S$ and the nutrient $N$.
We recall that $\sgn$ is the sign function.

\subsection{Theoretical bifurcation result}
\label{theobif}
Separately, the two subpopulations travel at different speeds and we name $\sigma_1$ the slow speed and $\sigma_2$ the fast one. 
The speeds $\sigma_i$, $i=1,2$, are computed thanks to the one species formula 
 \eqref{sigma_one_species} with the corresponding $\chi_i^S$ and $\chi_i^N$.
We denote $M_i$ the size of the subpopulation $i$, $\phi_{red}=\frac{M_2}{M_1+M_2}$ the fraction of the fast subpopulation (bacteria mCherry) and $I_i$ the interval $[\chi_i^N-\chi_i^S,\chi_i^N+\chi_i^S]$ for $i=1,2$.\\
We can prove our \textbf{main result}:
if the following assumption holds
\begin{equation}\label{hyp}
 I_1 \cap I_2 \neq \varnothing, \quad \text{and } \sigma_2, \chi^N_2-\chi_2^S \not \in I_1 \cap I_2,\quad \text{where  } I_i:=[\chi_i^N-\chi_i^S,\chi_i^N+\chi_i^S],
\end{equation}
Then, there exists $\phi_{red}^{*} \in ]0,1[$ such that 
\begin{itemize}
 \item for $\phi_{red} \leq \phi_{red}^{*}$, there exist traveling pulses. Moreover, the speed of the wave $\sigma$ is between $\sigma_1$ and $\sigma_2$ and satisfies
 \begin{equation}\label{speed}
  (\sigma-\chi^N_1)+\chi_1^S\frac{\sigma}{\sqrt{\sigma^2+4\alpha D_S}}+\frac{\phi_{red}}{1-\phi_{red}}H(\sigma)\left((\sigma-\chi^N_2)+\chi_2^S\frac{\sigma}{\sqrt{\sigma^2+4\alpha D_S}}\right)=0,
 \end{equation}
 where $H$ is defined in  \eqref{eq_H}.
\item for $\phi_{red} >\phi_{red}^{*}$, there do not exist single-speed traveling pulses 
\end{itemize}
We remark that the speed of the wave $\sigma$ is given by an implicit equation depending on the parameters of the model and the subpopulation sizes $M_i$. 
Note that in the case of a single population ($\phi_{red}=0$), we recover the single-species equation for $\sigma$ \eqref{sigma_one_species}. We also notice that with our system's parameters (see table \ref{parameters}) condition \eqref{hyp} is satisfied and thus such a $\phi_{red}^{*}$ exists in our case.

\subsection{Numerical method}

In order to provide comparisons between the solutions of the mathematical model for two species \eqref{two_species_model} and the experimental data, we perform numerical simulations.
Equation~\eqref{two_species_model} is discretized by a finite difference semi-implicit scheme. 
Such schemes are employed to solve numerically advection-diffusion equations.
This allows to avoid a too restrictive CFL condition imposed by diffusive terms. 
It consists in using an implicit time integration scheme for diffusive terms and explicit time integration for other terms.
Central finite differences are used for diffusive terms whereas a finite volume approach allows to discretize advection terms of the equation for $\rho_i$.

Let us consider a cartesian grid of space step $\Delta x$ and time step $\Delta t$, then 
we denote $x_k =k\Delta x$, for $k\in \Z$, and $t^n=n\Delta t$, for $n\in \N$.
For $k\in \Z$ and $n\in \N$, we consider approximation of $\rho_i(t^n,x_k)$ for $i=1,2$, 
$S(t^n,x_k)$, and $N(t^n,x_k)$ by $\rho_{i,k}^{n}$, $i=1,2$, $S_k^n$,
and $N_k^n$, respectively.

For a given $n\in \N$, assume that $S_k^n,N_k^n$ are known at time $t^n$ for all $k\in \Z$.
Then, $S_k^{n+1}$ and $N_k^{n+1}$ are computed thanks to the following iterative process
\begin{equation*}
\begin{aligned}
  \frac{S_k^{n+1}-S_k^n}{\Delta t}&=\frac{D_S}{{\Delta x}^2}\left(S_{k+1}^{n+1}-2S_{k}^{n+1}+S_{k-1}^{n+1}\right)-\alpha S_k^{n+1}+\rho_{1,k}^{n+1}+\rho_{2,k}^{n+1},\\
 \frac{N_k^{n+1}-N_k^n}{\Delta t}&=\frac{D_N}{{\Delta x}^2}\left(N_{k+1}^{n+1}-2N_{k}^{n+1}+N_{k-1}^{n+1}\right)-\gamma^1 N_k^{n+1}\rho_{1,k}^{n}-\gamma^2 N_k^{n+1}\rho_{2,k}^n .
\end{aligned}
\end{equation*}
This boils down to solving a linear system.

For the equation for $\rho_i$, $i=1,2$, we proceed differently. The diffusion term is treated as before using an implicit discretization, whereas the advection term is discretized thanks to a finite volume method of upwind kind.
The overall discretization of $\rho_i$ reads:
\begin{equation}\label{num_rho}
 \frac{\rho^{n+1}_{i,k}-\rho^{n}_{i,k}}{\Delta t}=\frac{D_{i}}{\Delta x^2} (\rho^{n+1}_{i,k+1}-2\rho^{n+1}_{i,k}+\rho^{n+1}_{i,k-1})+\frac{1}{\Delta x}(F^{n}_{i,k+1/2}-F^{n}_{i,k-1/2}), \qquad i=1,2,
\end{equation}
where $F^{n}_{i,k+1/2}$ is given by 
\begin{equation}\label{num_rho_2}
 F^{n}_{i,k+1/2}=(a^{n}_{i,k}[S]+a^{n}_{i,k}[N])^{+}\rho^{n}_{i,k} -(a^{n}_{i,k+1}[S]+a^{n}_{i,k+1}[N])^{-}\rho^{n}_{i,k+1},
\end{equation}
and $a^+=\max\{0,a\}$ and $a^-=\max\{0,-a\}$ denote, respectively, the positive and negative part 
of a real $a$. The discretized velocities $a^{n}_{i,k}[S]$ and $a^{n}_{i,k}[N]$ are given by
\begin{equation*}
  a^{n}_{i,k}[S]=\chi_S^i \sgn\left(\frac{S_{k+1}^n-S_k^n}{\Delta x}\right),\qquad
  a^{n}_{i,k}[N]=\chi_N^i \sgn\left(\frac{N_{k+1}^n-N_k^n}{\Delta x}\right).
\end{equation*}
Since the scheme \eqref{num_rho} can be written under the form 
$\rho^{n+1}_{i,k}=\rho^{n}_{i,k}+ \frac{\Delta t}{\Delta x}(J_{k+1/2}-J_{k-1/2})$,
we verify easily by summing over $k\in\Z$ that the total mass is conserved.
Finally, we observe that the velocity field is discontinuous and we mention that 
in the case without diffusion ($D_{i}=0$), bacteria profiles may concentrate strongly into Dirac deltas.
The convergence of scheme \eqref{num_rho}--\eqref{num_rho_2},
even in this singular case with discontinuous velocities, 
in the sense of measures, has been studied in \cite{CEJLNV}
(we refer to \cite{FJNV} for the one species case).

Fig \ref{two_regimes} displays a comparison between the experimental results in Fig \ref{two_regimes}A and Fig \ref{two_regimes}B,
and the numerical simulations obtained with the above scheme in Fig \ref{two_regimes}C and Fig \ref{two_regimes}D.
The insets on the top left of each figure depict the spatial concentration profiles of each population
of bacteria at the instant corresponding to the white dashed line in the 
kymograph. Fig \ref{two_regimes}A and Fig \ref{two_regimes}C correspond to a ratio $\phi_{red}=10\%$, whereas for
Fig \ref{two_regimes}B and Fig \ref{two_regimes}D we have $\phi_{red}=90\%$. 
Comparing Fig \ref{two_regimes}A and Fig \ref{two_regimes}C, we observe that for low values of the ratio $\phi_{red}$ 
the matching between experimental and numerical results is very good, 
confirming our theoretical result.
For Fig \ref{two_regimes}C and Fig \ref{two_regimes}D, $\phi_{red}$ is beyond the threshold value $\phi_{red}^*$ and thus we do not 
have the existence of a traveling wave. This is illustrated here by the fact that the 
total population may split into several branches. 
In this case we have a fast mode (1) and a slow mode (2) that we will see again in Fig \ref{bifurcation}.

\section*{Quantitative and qualitative discussion}
\setcounter{section}{3}
\setcounter{subsection}{0}

In this work we derive a macroscopic model \eqref{two_species_model} from microscopic assumptions on the run and tumble motion of individual bacteria (in the spirit of \cite{VCJS}). This model extends the one proposed in \cite{VCJS} for the one species case. The analytical study 
of the model enables us to determine the profiles of the traveling wave solutions and to show that they travel with a speed $\sigma$ given by \eqref{speed}. Moreover, we prove the existence of a critical proportion $\phi^\ast_{red}$ of the red subpopulation above which the 
theoretical single traveling wave solution no longer exists.

\begin{figure}[ht!]
\centering
 \includegraphics[scale=0.5]{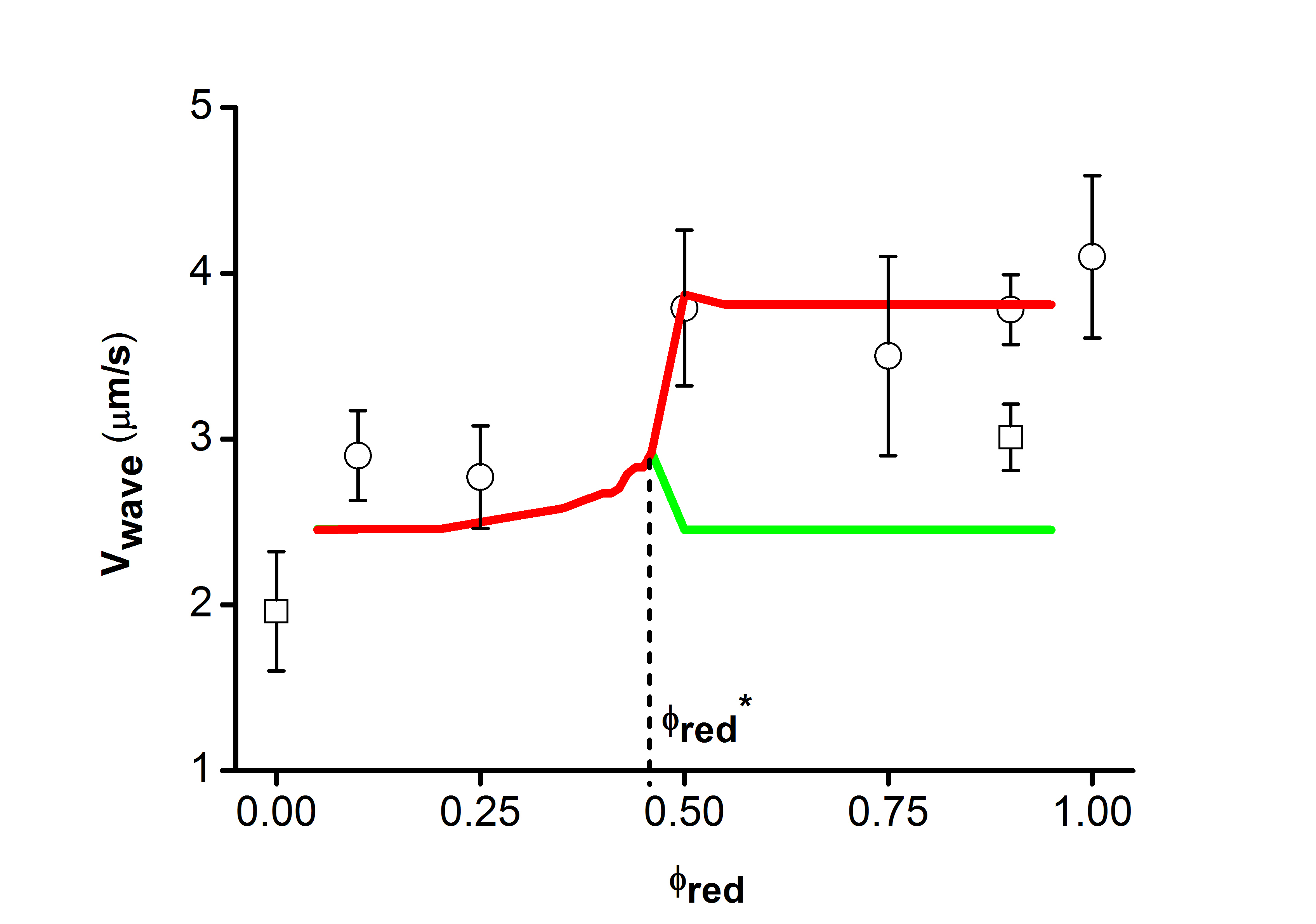}
 \caption{Mean wave velocity of red bacteria as a function of bacterial composition $\phi_{red}$:
the experimental results are represented by their means (open circles) and standard deviations (error bars). We use at least ten values for each point. 
These experimental results are compared to the simulation (red and green curves). 
One observes a separation of the two strains around $\phi_{red}=50\%$. The square corresponds to the mean velocity of the green wave observed in Fig 3D.}
 \label{bifurcation}
\end{figure}

Beyond this value $\phi^\ast_{red}$, we clearly observe in Fig \ref{two_regimes} that subpopulations 1 and 2 split.
After a transitory regime, they move separately at their own speed. We notice that it might take a long time for the separation to be completed in order to have a well defined speed for each pulse.
Then, to define these speeds we consider the density profiles of the green and red subpopulations at each later time and notice that both subpopulations have a clear peak around their maximum density. 
It is the spatial position of this maximum point at each time that is used to define the position of the corresponding subpopulation pulse. This position is then used to compute the speed of each 
subpopulation pulse.
These speeds are reported in Fig \ref{bifurcation} where a comparison between the 
numerical results and the experimental results is provided for different values of $\phi_{red}$.

The numerical simulations (Fig \ref{bifurcation}) show that, for small values of $\phi_{red}$, the two subpopulations travel together in a single wave of speed $\sigma$ given in Equation \eqref{speed} as predicted by our theoretical study. 
For higher values of $\phi_{red}$ there is clearly a bifurcation at the proportion $\phi^\ast_{red}$ beyond which 
the two subpopulations travel in separate waves with different speeds. The experimental results (Fig \ref{two_regimes}) confirm this separation of the two subpopulations for high values of $\phi_{red}$ but our experiments are not sufficiently precise
to enable us to look for the experimental bifurcation point (for intermediate values of $\phi_{red}$ the experimental data on the population density are diffuse and do not have clear peaks).

In Fig \ref{two_regimes}, we notice that for $\phi_{red}=0.9$ the experimental point corresponding to the red subpopulation is below the red curve and the one 
corresponding to the Green subpopulation is above the Green curve. This should be due to the fact that our experimental channel is probably not long enough. 
In fact, in the simulation (see Fig \ref{two_regimes}D), the experimental setting corresponds to the white rectangle in the lower left corner,
while the full separation of the pulses only happens for later times (further to the right in Fig \ref{two_regimes}D). Therefore, the experimental measurement might still be influenced 
by the transitory regime which would lead to overestimating the speed of the green pulse (due to the presence or interaction with the red subpopulation).

As observed in Fig \ref{two_regimes}., experimental results are not sufficiently precise 
for intermediate values of $\phi_{red}$, close to the bifurcation point.
In particular, we notice that for $\phi_{red}=75\%$ (and also for $\phi_{red}=50\%$), 
the error bar in the measurement of the speed of the front is larger than for other values.
It is due to the fact that for this intermediate value of $\phi_{red}$, our experimental results 
do not provide a clear separation into two branches of the wave front.

We also remark that the environment is changed after the passing of one wave of bacteria.
Consequently, the second wave of bacteria evolves in a different biochemical environment from the one seen by first wave.
This is not taken into account in our mathematical model.
We hope that it will be possible to do more detailed experiments (in particular having a better knowledge of the changes in the properties of the medium when it is crossed by the first wave) in the near future 
since they would provide us with important information 
to support the predictions of our model (and make it evolve if necessary).

As in \cite{VCJS}, we could do the analytical study only in the case where the functions $u_i[S],u_i[N]$ have the particular form \eqref{signals}. It would be extremely interesting, but very challenging, to be able to extend this study to more general cases.

\section*{Materials and Methods}
\setcounter{section}{4}
\setcounter{subsection}{0}

\subsection{Bacterial Strain and Cell culture}
We used the strain RP437 considered wild type for motility and chemotaxis. The strains were transformed by heat shock with PZE1R-GFP and PZE1R-mCherry plasmids. 
Cells were cultured in 3 mL LB medium (Sigma) with ampicillin at 33 $^{\circ}{\rm C}$, with shaking, up to mid-exponential phase (Optical density $OD_{600}=0.5$), and re-suspended after centrifugation in the medium used for the experiments: M9 Minimal Salts, 
5$\times$ supplemented with 1 g·L−1 Bacto™ Casamino Acids (both from Difco Laboratories, Sparks), 4 g·L−1 D-Glucose, and 1 mM MgSO4. 
The two types of bacteria were cultured independently before being mixed at the desired ratio at a final concentration corresponding to $OD_{600}=0.5$.

\subsection{Micro fabrication and centrifugation}
The micro-channels were prepared using usual soft lithography techniques \cite{Xia}. 100$\mu$m-high patterns were micro-fabricated on silicon wafers using SU-8 100 resin (MICROCHEM). 
The PDMS was molded on the wafer and peeled off after curing. A clean glass slide and the micro patterned PDMS were plasma treated for 30s and directly placed in contact thereby forming an array of 8 PDMS/glass parallel micro-channels (width=500$\mu$m, height=100$\mu$m, length=1.8cm). 
They were then filled by capillarity with the homogeneous suspension of motile bacteria and sealed with a fast curing epoxy resin.
The glass slide was gently centrifuged (800rpm, rotor diameter 20cm) at room temperature for half an hour. The bacteria accumulated at one end of the channels and stayed motile.

\subsection{Video Microscopy}
The channels were then immediately placed in a closed chamber maintained at constant temperature (33 $^{\circ}{\rm C}$). 
Few minutes after centrifugation had been stopped the concentration waves of bacteria started to propagate inside the channels (Fig \ref{one_species}). 
The observations were performed with a Leica MZ16F stereomicroscope equipped with two fluorescence sets: a green one, GFP2 (Leica) Ex480/Em510 and a red one, G (Leica) Ex546/Em590. 
Images were recorded on a CCD camera (CoolSnapHQ, Roper Scientific) at a frame rate of one image per minute (switching every minute from one fluorescence channel to the other). The image stacks were then post-processed using ImageJ and Matlab.

\subsection{Analytical forms of $\rho_i$ and $S$}
Experiments show that bacteria are concentrated locally in space while traveling. Therefore, traveling pulses (see \cite{Stevens,Nagai}) are particularly interesting to study.
By definition, we say that Equation \eqref{two_species_model} admits traveling pulses if and only if $\rho_i,S,N$ are traveling waves i.e functions satisfying the ansatz
\begin{equation*}
\rho_1(t,x)=\tilde \rho_1(z),\quad \rho_2(t,x)=\tilde \rho_2(z),\quad S(t,x)=\tilde S(z),\quad N(t,x)=\tilde N (z),\quad \text{where } z=x-\sigma t,
\end{equation*}
with $\sigma$ being the speed of the wave. The unknowns of the problem are $\sigma$ and the one-variable functions $\tilde S,\tilde N,\tilde \rho_1,\tilde \rho_2$, where $\tilde \rho_1,\tilde \rho_2$ are pulses as defined below.\medskip\\
\textbf{Definition}\\
A pulse is defined as a real-valued function which is increasing for negative values of $z$, decreasing for positive ones and decays to zero at infinity.\medskip\\
Plugging these expressions into Equation \eqref{two_species_model} and dropping the tilde over variables yields 
\begin{equation}
\begin{cases}
-\sigma {(\rho_1)}^{'}=D_1 {(\rho_1)}^{''}-\left ( \rho_1 \left (u_1[S]+ u_1[N] \right) \right )',\\
-\sigma {(\rho_2)}^{'}=D_2 {(\rho_2)}^{''}-\left( \rho_2 \left(u_2[S]+u_2[N]\right) \right )',\\
-\sigma {S}^{'}=D_{S} {S}^{''}-\alpha  S+\rho_1+\rho_2,\\
-\sigma {N}^{'}=D_{N} {N}^{''}-\gamma_1 \rho_1  N -\gamma_2 \rho_2 N.\\
\end{cases}
\label{eq_traveling_wave}
\end{equation}
Looking for a pulse, we prescribe the following boundary conditions
\begin{equation}\label{bc}
 \rho_1(\pm \infty)=0,\quad \rho_2(\pm \infty)=0,\quad S(\pm \infty)=0.
\end{equation}
The fact that $\rho_i$ are pulses implies that (see subsection \ref{secA1})
\begin{equation}\label{pulse}
 \begin{cases}
  & (\rho_i)^{'}(\pm \infty)=0,\quad S'(\pm \infty)=0,\\
  & \partial_z S >0 \quad \text{for } z>0,\quad \partial_z S<0 \quad \text{for } z>0,\\
  & \partial_z N >0 \quad \text{for } z \in \mathbb{R}.
 \end{cases}
\end{equation}
Therefore, $u_i[N]$ and $u_i[S]$ are given by 
\begin{equation*}
 u_i[S]=-\chi_i^S \sgn(z),\quad u_i[N]=\chi_i^N \quad \text{for }i=1,2.
\end{equation*}
We integrate equations for $\rho_i$ in \eqref{eq_traveling_wave} and use \eqref{pulse} to obtain 
\begin{equation}\label{eq_profile}
\begin{cases}
D_1{\rho_1}'=\left(u_1[S]+u_1[N]-\sigma \right) \rho_1,\\
D_2{\rho_2}'=\left(u_2[S]+u_2[N]-\sigma\right )\rho_2.
\end{cases}
\end{equation} 
We deduce the analytical forms of $\rho_i$ for $i=1,2$ and $S$.

\begin{equation}\label{eq_rho}
{\rho}_i=
\left 
\{
\begin{aligned}
\rho_i^M \exp(\lambda_i^{-}z),\quad \lambda_i^{-}=\frac{\chi_i^N+\chi_i^S-\sigma}{D_i}>0, \quad & \text{for } z<0, \\
\rho_i^M \exp(\lambda_i^{+}z),\quad \lambda_i^{+}=\frac{\chi_i^N-\chi_i^S-\sigma}{D_i}<0, \quad & \text{for } z>0, 
\end{aligned}\quad i=1,2.
\right.
\end{equation}
From the equation satisfied by S, we deduce that
\begin{equation}\label{eq_S}
 S(z)=\int_{-\infty}^{+\infty} K(z-y)\left(\rho_1(y)+\rho_2(y)\right)dy,
\end{equation}
where $K$ is given by 
\begin{equation}\label{eq_K}
 K=\exp\left(-\frac{\sigma}{2D_S}z-\frac{\sqrt{\sigma^2+4\alpha D_S}}{2D_S}\abs{z}\right).
\end{equation}
\subsection{Speed of the wave $\sigma$}

Since S is maximal at $z=0$, by definition of a pulse, we should have $S'(0)=0$.
Differentiating $S$ in \eqref{eq_S} gives
\begin{equation*}
 S'(0)=\int_{-\infty}^{+\infty} K'(x)(\rho_1(-x)+\rho_2(-x))dx.
\end{equation*}
We split this integral into two parts $S'_{-}$ and $S'_{+}$:
\begin{equation*}
\begin{aligned}
 S'(0)=&\int_{-\infty}^{0} K'(x)(\rho_1(-x)+\rho_2(-x))dx+\int_{0}^{+\infty} K'(x)(\rho_1(-x)+\rho_2(-x))dx.\\
 S'(0)=& S'_{-}+S'_{+}.
\end{aligned}
\end{equation*}
From subsection \ref{secA21}, $S'_{-}$ and $S'_{+}$ are given by 
\begin{equation*}
\begin{aligned}
&S'_{-}=\rho_1^M \frac{-\sigma+\sqrt{\sigma^2+4\alpha D_S}}{-\sigma+\sqrt{\sigma^2+4\alpha D_S}-2D_S\lambda_1^{+}}+\rho_2^M \frac{-\sigma+\sqrt{\sigma^2+4\alpha D_S}}{-\sigma+\sqrt{\sigma^2+4\alpha D_S}-2D_S\lambda_2^{+}},\\
&S'_{+}=\rho_1^M\frac{\sigma+\sqrt{\sigma^2+4\alpha D_S}}{-\sigma-\sqrt{\sigma^2+4\alpha D_S}-2D_S\lambda_1^{-}}+\rho_2^M\frac{\sigma+\sqrt{\sigma^2+4\alpha D_S}}{-\sigma-\sqrt{\sigma^2+4\alpha D_S}-2D_S\lambda_2^{-}}.
\end{aligned}
\end{equation*}
Putting together $S'_{-}$ and $S'_{+}$ gives 
\begin{equation*}
 \begin{aligned}
  S'(0)=&\rho_1^M \frac{c_1}{\left(-\sigma-\sqrt{\sigma^2+4\alpha D_S}-2D_S\lambda_1^{-} \right) \left(-\sigma+\sqrt{\sigma^2+4\alpha D_S}-2D_S\lambda_1^{+}\right)}\\
        +&\rho_2^M \frac{c_2}{\left(-\sigma-\sqrt{\sigma^2+4\alpha D_S}-2D_S\lambda_2^{-}\right) \left(-\sigma+\sqrt{\sigma^2+4\alpha D_S}-2D_S\lambda_2^{+}\right)},
 \end{aligned}
\end{equation*}
with $c_i$ given by
\begin{equation}\label{eq_c}
c_i=4\frac{D_S}{D_i}\left (\chi_i^S \sigma+ \left(\sigma-\chi_i^N\right)\sqrt{\sigma^2+4\alpha D_S}\right),\quad i=1,2.
\end{equation}
Due to \eqref{pulse}, $S$ is maximal for $z=0$, then $S'$ vanishes at $0$ and we obtain the equation for $\sigma$
\begin{multline}\label{eq_sigma_1}
\left(\chi_1^S \sigma+ (\sigma-\chi^N_1)\sqrt{\sigma^2+4\alpha D_S}\right)+\left(\chi_2^S \sigma+ (\sigma-\chi^N_2)\sqrt{\sigma^2+4\alpha D_S}\right) \times\\
\frac{D_1\rho_2^M}{D_2 \rho_1^M}\frac{\left(-\sigma-\sqrt{\sigma^2+4\alpha D_S}-2D_S\lambda_1^{-}\right) \left(-\sigma+\sqrt{\sigma^2+4\alpha D_S}-2D_S\lambda_1^{+}\right)}{\left(-\sigma-\sqrt{\sigma^2+4\alpha D_S}-2D_S\lambda_2^{-}\right) \left(-\sigma+\sqrt{\sigma^2+4\alpha D_S}-2D_S\lambda_2^{+}\right)}=0.
\end{multline}
We recall that $\lambda^i_{\pm}$ are given by \eqref{eq_rho}. Unknowns $\rho_i^M$ are obtained thanks to the conservation of the total subpopulation $M_i$.
Indeed, from the conservative form of equations for $\rho_i$ \eqref{two_species_model}, it follows that for all $t\geq 0$
\begin{equation*}
 M_i=\int_{-\infty}^{+\infty}\rho_i(x,t)dx=\int_{-\infty}^{+\infty} \rho_i^{ini}(x)dx,\quad i=1,2,
\end{equation*}
with $\rho_i^{ini}=\rho_i(x,t=0)$ the initial profile of $\rho_i$. We deduce
%This leads to 
\begin{equation*}
 M_i=\int_{-\infty}^{+\infty} \rho_i(z)dz=\rho_i^M \frac{\chi_i^S D_i}{(\chi_i^S)^2-(\sigma-\chi_i^N)^2}.
\end{equation*}
Replacing $\lambda_i^{+}$, $\lambda_i^{-}$ by their values in \eqref{eq_rho} and using the previous relationship, \eqref{eq_sigma_1} becomes
\begin{equation}\label{eq_sigma2}
(\sigma-\chi^N_1)+\chi_1^S\frac{\sigma}{\sqrt{\sigma^2+4\alpha D_S}}+\frac{\phi_{red}}{1-\phi_{red}}H(\sigma)\left((\sigma-\chi^N_2)+\chi_2^S\frac{\sigma}{\sqrt{\sigma^2+4\alpha D_S}}\right)=0,
\end{equation}
with 
\begin{equation}\label{eq_H}
\begin{aligned}
H(\sigma)&=\frac{\chi_1^S D_1}{\chi_2^S D_2}\frac{(\chi_2^S)^2-(\sigma-\chi^N_2)^2}{(\chi_1^S)^2-(\sigma-\chi^N_1)^2}\frac{h_1(\sigma)}{h_2(\sigma)},\\
h_i(\sigma)&= \sigma^2(\frac{D_S}{D_i}-1)+(1-2\frac{D_S}{D_i})\sigma \chi_i^N-\chi_i^S\sqrt{\sigma^2+4\alpha D_S}+\frac{D_S}{D_i}((\chi_i^N)^2-(\chi_i^S)^2)-\alpha D_i, \quad i=1,2.
\end{aligned}
\end{equation}
From subsection \ref{secA3}, we have that $\sigma$ belongs to $\Omega=I_1 \cap I_2 \cap (\sigma_1,\sigma_2)$.
Since $\sigma_2$ and $\chi^N_2-\chi_2^S$ do not belong to $\Omega$, then \eqref{eq_sigma2} rewrites 
\begin{equation*}
 G(\sigma)=\frac{\phi_{red}}{1-\phi_{red}},
\end{equation*}
where $G$ is a positive function bounded over $\Omega$ given by 
\begin{equation*}
 G(\sigma)=-\frac{\chi_1^S \sigma+ (\sigma-\chi^N_1)\sqrt{\sigma^2+4\alpha D_S}}{H(\sigma)\left(\chi_2^S \sigma+ (\sigma-\chi^N_2)\sqrt{\sigma^2+4\alpha D_S}\right)}.
\end{equation*}
The function $G$ admits a maximum which is finite. Therefore, there exists $\phi_{red}^{*} \in (0,1)$ such that 
\begin{equation*}
\max_{\Omega}{G(\sigma)}=\frac{\phi_{red}^*}{1-\phi_{red}^*}.
\end{equation*}
For more details, we refer to subsection \ref{secA3}.

\subsection{Derivation of the two-species macroscopic model}\label{sec:kin}
In this subsection, we derive formally macroscopic equations \eqref{two_species_model} from the kinetic descriptions of individual motion of bacteria \cite{OthmerHill,CMPS,HKS}. This motion is a succession of run and tumble phases as observed in \cite{Berg}. 
During the run phase, bacteria move in straight lines and change their direction during the tumble phase. The Othmer-Dunbar-Alt model (see \cite{OthmerAlt,Wolfgang,Saragosti}) gives the mathematical description of the individual behavior.
It describes the dynamics of the distribution density $f_i(x,v,t)$ of cells at position $x \in \mathbb{R}^d$ at time $t$ having speed $v \in V$, where $V$ is a bounded, symmetric and rotationally invariant domain of $\mathbb{R}^d$.
It reads
\begin{equation}\label{kinetics}
 \begin{aligned}
  & \partial_t f_i  +v\cdot \nabla_x f_i =\frac{1}{2}\int_{V}\left(T_i[S](x,v,v',t)f_i(x,v',t)-T_i[S](x,v',v,t)f_i(x,v,t)\right)dv'\\
  & \phantom{\partial_t f_i  +v\cdot \nabla_x f_i}+\frac{1}{2}\int_{V}\left(T_i[N](x,v,v',t)f_i(x,v',t)-T_i[N](x,v',v,t)f_i(x,v,t)\right)dv',\quad \text{for }i=1,2,\\
 \end{aligned}
\end{equation}
where $T_i[S](x,v,v',t)$, $T_i[N](x,v,v',t)$ stand for the amount of bacteria reorienting from the direction $v'$ to $v$ under the influence of the chemoattractant and the nutrient.
They are usually called tumbling kernels.

As proposed in \cite{DolakSch}, we consider that bacteria have 
a small memory effect allowing them to sense the chemical concentrations along 
their trajectory and therefore to respond to gradients of concentrations.
Then the tumbling kernels $T_i[S]$ and $T_i[N]$ are given by
\begin{equation*}
\left\{
\begin{aligned}
T_i[S]&=\Phi_i[S](\partial_t S+v'\cdot \nabla_x S),\\
T_i[N]&=\Phi_i[N](\partial_t N+v'\cdot \nabla_x N),
 \end{aligned}
 \right.
\end{equation*}
Let us introduce $\varepsilon$ the mean free path, i.e. the average distance between two successive 
tumblings, usually $\varepsilon\ll 1$. When the taxis (e.g. chemotaxis) is small compared to the unbiased movement of cells, 
we perform a diffusive scaling ($\tilde{x}=\varepsilon x,\tilde{t}=\varepsilon^2 t$) and get
\begin{equation*}
 \begin{aligned}
  & \partial_t f_i^{\varepsilon}  +\frac{v}{\varepsilon}\cdot \nabla_x f_i^{\varepsilon} =\frac{1}{2\varepsilon^2} \int_{V}\Phi_i^{\varepsilon}[S](\partial_t S+v'\cdot \nabla_x S)f_i^{\varepsilon}(x,v',t)dv'-\abs{V}\Phi_i^{\varepsilon}[S](\partial_t S+v\cdot \nabla_x S)f_i^{\varepsilon}(x,v,t),\\
  & \phantom{\partial_t f_i^{\varepsilon}  +\frac{1}{\varepsilon} v\cdot \nabla_x f_i^{\varepsilon}}+\frac{1}{2\varepsilon^2} \int_{V}\Phi_i^{\varepsilon}[N](\partial_t N+v'\cdot \nabla_x N)f_i^{\varepsilon}(x,v',t)dv'-\abs{V}\Phi_i^{\varepsilon}[N](\partial_t N+v\cdot \nabla_x N)f_i^{\varepsilon}(x,v,t).\\
 \end{aligned}
\end{equation*}
Since tumbling kernels are perturbations of constant tumbling rates for E.Coli, we can assume that 
\begin{equation*}
\left\{
 \begin{aligned}
\Phi_i^{\varepsilon}[S](\partial_t S+v\cdot\nabla_x S)=&\psi_i \left(1+\varepsilon \phi_i^S(\varepsilon \partial_t S+v\cdot \nabla_x S)\right),\\
\Phi_i^{\varepsilon}[N](\partial_t N+v\cdot\nabla_x N)=&\psi_i \left(1+\varepsilon \phi_i^N(\varepsilon \partial_t N+v\cdot \nabla_x N)\right),
 \end{aligned}
 \right.
\end{equation*}
where $\phi_i^S,\phi_i^N$ are decreasing functions.\\
We write expansions of $f_i^{\varepsilon}$ when $\varepsilon$ tends to zero.
\begin{equation*}
 f_i^{\varepsilon}=f_i^0+\varepsilon f_i^1,\quad \text{for }i=1,2.
\end{equation*}
Plugging these expansions in the equation for $f_i^{\varepsilon}$ gives at the order $1/\varepsilon^2$
\begin{equation*}
 f_i^0(x,v,t)=\frac{\int_V f_i^0(x,v',t)dv'}{\abs{V}}=\frac{\rho_i^0(x,t)}{\abs{V}},\quad \text{for }i=1,2,
\end{equation*}
where $|V|$ is the measure of $V$.
At the order $1/\varepsilon$, we get 
\begin{multline*}
 f_i^1=\frac{1}{\abs{V}}\int_V f_i^1(x,v',t)dv'+\frac{1}{2\abs{V}}\int_V \left(\phi_i^S(v\cdot \nabla_x S^{0})+\phi_i^N(v\cdot \nabla_x N^{0})\right)f_i^0(x,v',t)dv'\\
 -\frac{1}{2}\left(\phi_i^S(v\cdot \nabla_x S^{0})+\phi_i^N(v\cdot \nabla_x N^{0}) \right)f_i^0(x,v,t)-\frac{v}{\abs{V}\psi_i}\cdot \nabla_x f_i^0,
\end{multline*}
where $S^{0}$ and $N^{0}$ are leading order terms of asymptotic expansions of respectively $S$ and $N$ and are solutions to equations \eqref{two_species_model} for $S$ and $N$ with $\rho_i=\rho_i^0$.\\
Integrating the equation for $f_i^{\varepsilon}$ over $V$ yields the following conservation equation for $\rho_i^{\varepsilon}=\int_V f_i^{\varepsilon}dv$
\begin{equation*}
 \partial_t \int_V f_i^{\varepsilon}dv+\frac{1}{\varepsilon}\nabla \cdot \left(\int_V v f_i^{\varepsilon}dv \right)=0.
\end{equation*}
From the asymptotic analysis carried out before, we have that 
\begin{equation*}
 \int_V f_i^{\varepsilon}(x,v,t)dv \rightarrow \rho_i^0(x,t).
\end{equation*}
Since $V$ is a symmetric bounded domain, we have the convergence of the scaled first moment 
\begin{equation*}
 \frac{1}{\varepsilon} \int_V v f_i^{\varepsilon}(x,v,t)dv \rightarrow \int_V v f_i^1(x,v,t)
\end{equation*}
with 
\begin{multline*}
 \int_V v f_i^1(x,v,t)=-\left(\int_V\frac{v}{2\abs{V}}\phi_i^S(v\cdot \nabla_x S^{0})dv\right) \rho_i^0(x,t) -\left(\int_V\frac{v}{2\abs{V}}\phi_i^N(v\cdot \nabla_x N^{0})dv\right) \rho_i^0(x,t)\\
 -\int_V\frac{v\otimes v}{\abs{V}^2\psi_i}\nabla_x \rho_i^0(x,t).
\end{multline*}
We finally obtain the equation for $\rho_i^0$ in \eqref{two_species_model}. This formal computation has been rigorously established by part of the authors in \cite{AEV}.

\subsection{Parameter estimation}

In this part, we discuss the estimation of parameters in the model \eqref{two_species_model}. We use a similar approach to the one-species model \eqref{one_species_model} one developped in \cite{VCJS}. Diffusion coefficients $D_S,D_1,D_2$ are measured experimentally (see \cite{Berg,Salman,Berg1}) whereas the other parameters are fitted thanks to experimental data.\\
The fact that the two-species speed formula \eqref{speed} extends the single-species one \eqref{sigma_one_species} leads us to fit separately parameters $\chi_1^N,\chi_1^S$ and 
$\chi_2^N,\chi_2^S$.\\ 
For the fitting of the couple of parameters $(\chi_i^S,\chi_i^N)$ for $i=1,2$, we use experimental data on the migration of the pure population $i$ which 
corresponds to extreme cases ($\phi_{red}=0$ and $\phi_{red}=100 \%$) and carry out the method used for the single-species case in \cite{VCJS}.\\
The double asymmetric exponential profile of bacteria measured by $\lambda_i^{\pm} \ (\lambda_i^{-}>0,\lambda_i^{+}<0)$ and the speed $\sigma_i$ allow us to obtain $\chi_i^S,\chi_i^N,\alpha$ as follows
\begin{equation*}
\begin{aligned}
\chi_i^S&=\frac{D_i\left(\lambda_i^{+}-\lambda_i^{-}\right)}{2},\\
\chi_i^N&=\sigma_i+\frac{D_i\left(\lambda_i^{+}+\lambda_i^{-}\right)}{2},\\
\alpha&=(\sigma_i)^2 \frac{-\lambda_i^{+}\lambda_i^{-}}{D_S(\lambda_i^{+}+\lambda_i^{-})^2}.\\
\end{aligned}
\end{equation*}
Parameters $D_S,D_i,\chi_i^S,\chi_i^N,\alpha$ are found in Table \ref{parameters}.

\section*{Technical results}
\setcounter{section}{5}
\setcounter{subsection}{0}

In this section, we prove several technical results. It is written to provide the details of the proofs and can be skipped by less mathematically oriented readers.
\subsection{Proof of the result \eqref{pulse} on the signs of $\partial_z S$ and $\partial_z N$}
\label{secA1}
We prove that : if $\rho_1$ and $\rho_2$ are pulses, then $S$ is also a pulse and we have
\begin{equation*}
\left\{
\begin{aligned}
& \partial_z S>0,\quad \text{for }z<0,\\
& \partial_z S<0,\quad \text{for }z>0.
\end{aligned}
\right.
\end{equation*}
and 
\begin{equation*}
\partial_z N>0, \quad \text{for }z \in \mathbb{R}.
\end{equation*}
Let us prove this result. 
Equation \eqref{eq_S} says that $S=K*(\rho_1+\rho_2)$. Since $\rho_i$ decays to zero at infinity, we conclude that $S$ also decays to zero at infinity.
The equation for $S$ in \eqref{two_species_model} also implies that $S'$ admits a limit at infinity.
This limit has to be zero otherwise $S$ would not have a limit at infinity.
Since limits of $S$ are equal to zero at infinity, $S'$ vanishes at least once. Without loss of generality, let us suppose that $S'$ vanishes at $z=0$ 
if not, we can do a translation in $z$.
Differentiating the equation for $S$ in \eqref{two_species_model} yields
\begin{equation*}
-\sigma \partial_z S'=D_S \partial_{zz}(S')-\alpha S'+\partial_z \rho_1+\partial_z \rho_2
\end{equation*}
Since $\rho_i$ is a pulse, we have that $\partial_z \rho_i$ is positive in $(-\infty,0)$.
We get that 
\begin{equation*}
-\sigma \partial_z S'-D_S \partial_{zz} S'+\alpha S'\geq 0 \quad \text{in } (-\infty,0).
\end{equation*}
Consider $u=-S'$, then $u$ satisfies
\begin{equation*}
-\sigma \partial_z u-D_S \partial_{zz}u+\alpha u\leq 0,\quad (-\infty,0).
\end{equation*}
Multiplying this equation by $u^+:=\max(u,0)$ and integrating by parts yields
\begin{equation*}
-\sigma \int_{-\infty}^0 \partial_z u^+ u^+ +D_S \int_{-\infty}^0 {\partial_z u^{+}}^2+\alpha \int_{-\infty}^0 {u^{+}}^2 \leq 0.
\end{equation*}
Since the first term is zero 
\begin{equation*}
\int_{-\infty}^0 \partial_z u^+ u^+=[(u^+)^2]_{-\infty}^0=0,
\end{equation*}
we get 
\begin{equation*}
D_S \int_{-\infty}^0 {(\partial_z u^{+})}^2+\alpha \int_{-\infty}^0 {(u^{+})}^2 \leq 0.
\end{equation*}
Since we have the sum of two nonnegative terms, this implies that
\begin{equation*}
\int_{-\infty}^0 {(u^{+})}^2 \leq 0.
\end{equation*}
We have proved that $u^+=0$ for $z\in (-\infty,0)$ which by definition d of $u$ means that $\partial_x S>0$.\\
By a similar argument, we show that $\partial_z S<0$ for $z\in (0,\infty)$.\\
Now, we prove that $\partial_z N>0$.
By using the same technique, we can show that $\rho_i,i=1,2$ and $N$ are positive.
Denoting $u=-\partial_z N$ and using the positivity of $\rho_i$ and $N$ gives that 
\begin{equation*}
-\sigma u-D_N u'\geq 0,\quad \text{in }(-\infty,\infty).
\end{equation*}
By multiplying by $u^+$ and integrating by parts, one has
\begin{equation*}
\sigma \int_{-\infty}^{+\infty}(u^+)^2 \leq 0 \quad \text{which implies that }\partial_z N>0,\quad \text{in } (-\infty,-\infty).
\end{equation*}
\subsection{Detailed computation of $S'(0)$}
\label{secA2}
In this part, we provide missing computations of $S'_{-},S'_{+},c_i,h_i$ which are used to derive the dispersion relation \eqref{speed} in the Results section.
\subsubsection{Computations of $S'_{-},S'_{+}$}\label{secA21}
From the definition of $K$ \eqref{eq_K}, $K'$ reads
\begin{equation*}
K'=
\left\{
\begin{aligned}
& \left(-\frac{\sigma}{2D_S}+\frac{\sqrt{\sigma^2+4D_S \alpha}}{2D_S}\right) \exp \left(-\frac{\sigma}{2D_S}z+\frac{\sqrt{\sigma^2+4D_S \alpha}}{2D_S}z\right),\quad \text{for }z<0,\\
& \left(-\frac{\sigma}{2D_S}-\frac{\sqrt{\sigma^2+4D_S \alpha}}{2D_S}\right) \exp \left(-\frac{\sigma}{2D_S}z-\frac{\sqrt{\sigma^2+4D_S \alpha}}{2D_S}z\right),\quad \text{for }z>0.\\
\end{aligned}
\right.
\end{equation*}
We have 
\begin{equation*}
\int_{-\infty}^0 K'(x)\rho_i(-x)dx=\rho_i^M\frac{-\sigma \sqrt{\sigma^2+4\alpha D_S}}{2D_S}\int_{-\infty}^0\exp \left(-\frac{\sigma}{2D_S}x+\frac{\sqrt{\sigma^2+4D_S \alpha}}{2D_S}x-\lambda_i^{+} x \right) dx.
\end{equation*}
It follows that 
\begin{equation*}
\int_{-\infty}^0 K'(x)\rho_i(-x)dx=\rho_i^M\frac{\frac{-\sigma+\sqrt{\sigma^2+4\alpha D_S}}{2D_S}}{\frac{-\sigma+\sqrt{\sigma^2+4\alpha D_S}}{2D_S}-\lambda_i^{+}} \times \left[\exp\left(-\frac{\sigma}{2D_S}x+\frac{\sqrt{\sigma^2+4D_S \alpha}}{2D_S}x-\lambda_i^{+} x\right)\right]_{-\infty}^0.
\end{equation*}
We conclude that 
\begin{equation*}
\int_{-\infty}^0 K'(x)\rho_i(-x)dx=\rho_i^M\frac{-\sigma+\sqrt{\sigma^2+4\alpha D_S}}{-\sigma+\sqrt{\sigma^2+4\alpha D_S}-2D_S\lambda_i^{+}}.
\end{equation*}
In the same way, we can prove that
\begin{equation*}
\int_0^{+\infty} K'(x)\rho_i(-x)dx=\rho_i^M\frac{\sigma+\sqrt{\sigma^2+4\alpha D_S}}{-\sigma-\sqrt{\sigma^2+4\alpha D_S}-2D_S\lambda_i^{-}}.
\end{equation*}
Finally, $S'_{-}$ and $S'_{+}$ read 
\begin{equation*}
\begin{aligned}
 S'_{-}&=\int_{-\infty}^{0} K'(x)(\rho_1(-x)+\rho_2(-x))dx,\\
       &= \rho_1^M \frac{-\sigma+\sqrt{\sigma^2+4\alpha D_S}}{-\sigma+\sqrt{\sigma^2+4\alpha D_S}-2D_S\lambda_1^{+}}+\rho_2^M \frac{-\sigma+\sqrt{\sigma^2+4\alpha D_S}}{-\sigma+\sqrt{\sigma^2+4\alpha D_S}-2D_S\lambda_2^{+}},
\end{aligned}
\end{equation*}
\begin{equation*}
\begin{aligned}
 S'_{+}&=\int_{0}^{\infty} K'(x)(\rho_1(-x)+\rho_2(-x))dx\\
       &=\rho_1^M\frac{\sigma+\sqrt{\sigma^2+4\alpha D_S}}{-\sigma-\sqrt{\sigma^2+4\alpha D_S}-2D_S\lambda_1^{-}}+\rho_2^M\frac{\sigma+\sqrt{\sigma^2+4\alpha D_S}}{-\sigma-\sqrt{\sigma^2+4\alpha D_S}-2D_S\lambda_2^{-}}.
\end{aligned} 
\end{equation*}
\subsubsection{Computations of $c_i$ and $h_i$}\label{secA22}
Coefficients $c_i$ are defined in \eqref{eq_c} 
\begin{multline*}
c_i=\left(-\sigma+\sqrt{\sigma_2+4 \alpha D_S}\right) \left(-\sigma-\sqrt{\sigma^2+4\alpha D_S}-2D_S \lambda_i^{-}\right)\\
+\left(\sigma+\sqrt{\sigma^2+4\alpha D_S} \right)\left(-\sigma +\sqrt{\sigma^2+4\alpha D_S}-2D_S \lambda_i^{+}\right).
\end{multline*}
Expanding $c_i$ yields
\begin{multline*}
c_i=(-\sigma+\sqrt{\sigma^2+4\alpha D_S})(-\sigma-\sqrt{\sigma^2+4\alpha D_S})+(\sigma+\sqrt{\sigma^2+4\alpha D_S})(-\sigma+\sqrt{\sigma^2+4\alpha D_S})\\
+2D_S(\sigma-\sqrt{\sigma^2+4\alpha D_S})\lambda_i^{-}-2D_S(\sigma+\sqrt{\sigma^2+4\alpha D_S})\lambda_i^{+}.
\end{multline*}
Recalling expressions of $\lambda_i^{\pm}$
\begin{equation*}
\left\{
\begin{aligned}
\lambda_i^{-}&=\frac{\chi_i^N+\chi_i^S-\sigma}{D_i},\\
\lambda_i^{+}&=\frac{\chi_i^N-\chi_i^S-\sigma}{D_i},\\
\end{aligned}
\right.
\end{equation*}
we can write
\begin{equation*}
c_i=2D_S(\sigma-\sqrt{\sigma^2+4\alpha D_S}) \left\{ \frac{\chi_i^N-\sigma}{D_i}+\frac{\chi_i^S}{D_i}\right\}-2D_S(\sigma+\sqrt{\sigma^2+4\alpha D_S}) \left\{ \frac{\chi_i^N-\sigma}{D_i}-\frac{\chi_i^S}{D_i}\right\}.
\end{equation*}
After simplification, we obtain
\begin{equation*}
c_i=-4\frac{D_S}{D_i}(\chi_i^N-\sigma)\sqrt{\sigma^2+4\alpha D_S}+4\frac{D_S}{D_i}\chi_i^S \sigma.
\end{equation*}
We now compute the functions $h_i$ defined by:
\begin{equation*}
h_i(\sigma):=\left(-\sigma-\sqrt{\sigma^2+4\alpha D_S}-2D_S\lambda_i^{-}\right)\left(-\sigma+\sqrt{\sigma^2+4\alpha D_S}-2D_S \lambda_i^{+}\right).
\end{equation*}
The expansion of $h_i$ leads to
\begin{multline*}
h_i(\sigma)=\left(-\sigma-\sqrt{\sigma^2+4\alpha D_S}\right)\left(-\sigma+\sqrt{\sigma^2+4\alpha D_S}\right)\\
+2D_S\lambda_i^{+}\left(\sigma+\sqrt{\sigma^2+4\alpha D_S}\right)+2D_S\lambda_i^{-}\left(\sigma-\sqrt{\sigma^2+4\alpha D_S}\right)+4D_S^2\lambda_i^{+} \lambda_i^{-}.
\end{multline*}
Straightforward computations give
\begin{equation*}
h_i(\sigma)=-4\alpha D_S+4\frac{D_S}{D_i}\sigma(\chi_i^N-\sigma)-4\frac{D_S}{D_i}\chi_i^S\sqrt{\sigma^2+4\alpha D_S}+4(\frac{D_S}{D_{\rho}^i})^2\left((\chi_i^N-\sigma)^2-(\chi_i^S)^2\right).
\end{equation*}
We obtain after further simplifications
\begin{equation*}
h_i(\sigma)=4\frac{D_S}{D_i}\left(\sigma^2(\frac{D_S}{D_i}-1)+\chi_i^N \sigma(1-2\frac{D_S}{D_i})-\chi_i^S\sqrt{\sigma^2+4\alpha D_S}+\frac{D_S}{D_i}\left((\chi_i^N)^2-(\chi_i^S)^2\right)-\alpha D_i\right).
\end{equation*}
\subsection{Complete analysis of traveling pulses}\label{secA3}
In this part, we detail the study of the dispersion relation presented briefly in subsection \ref{theobif}.
We suppose the existence of traveling pulses, which implies the double exponential shapes of $\rho_i$ as given in \eqref{eq_rho} with
\begin{equation*}
\left\{
\begin{aligned}
\lambda_i^{-}&=\frac{\chi_i^N-\sigma+\chi_i^S}{D_{\rho}^i}>0\\
\lambda_i^{+}&=\frac{\chi_i^N-\sigma-\chi_i^S}{D_{\rho}^i}<0
\end{aligned}
\right.
\Rightarrow \chi_i^N-\chi_i^S<\sigma<\chi_i^N+\chi_i^S
\end{equation*}
The non-emptiness of the intersection of $I_1$ and $I_2$ in assumption \eqref{hyp} is the first condition to have the existence of traveling pulses.
This gives the interval to which $\sigma$ must belong. The dispersion relation will be studied in this interval.
For $\sigma$ in this interval, $\lambda_i^{-}$ is positive and $\lambda_i^{+}$ negative, then
\begin{equation*}
h_i(\sigma)=\left(\underbrace{-\sigma-\sqrt{\sigma^2+4\alpha D_S}}_{<0}-2D_S\lambda_i^{-}\right)\left(\underbrace{-\sigma+\sqrt{\sigma^2+4\alpha D_S}}_{>0}-2D_S\lambda_i^{+}\right) <0.
\end{equation*}
Moreover, 
\begin{equation*}
\chi_i^N-\sigma \in (-\chi_i^S,\chi_i^S) \Rightarrow (\chi_i^S)^2-(\chi_i^N-\sigma)^2>0.
\end{equation*}
Therefore,
\begin{equation*}
H(\sigma)=\frac{\chi_1^S D_1}{\chi_2^S D_2}\frac{(\chi_2^S)^2-(\chi^N_2-\sigma)^2}{(\chi_1^S)^2-(\chi^N_1-\sigma)^2}\frac{h_1(\sigma)}{h_2(\sigma)}>0.
\end{equation*}
We recall the dispersion relation \eqref{speed} obtained in the section Results,
\begin{equation*}
\left((\sigma-\chi^N_1)+\chi_1^S\frac{\sigma}{\sqrt{\sigma^2+4\alpha D_S}}\right)+\frac{\phi_{red}}{1-\phi_{red}}H(\sigma)\left((\sigma-\chi^N_2)+\chi_2^S\frac{\sigma}{\sqrt{\sigma^2+4\alpha D_S}}\right)=0.
\end{equation*}
Denote $g_i$ the map $\sigma \longmapsto (\sigma-\chi_i^N)+\chi_i^S\frac{\sigma}{\sqrt{\sigma^2+4\alpha D_S}}$ defined over $I_i$.
Then $g_i$ is an increasing map
\begin{equation*}
(g_i)^{'}(\sigma)=1+\frac{4\alpha D_S}{(\sigma^2+4\alpha D_S)\sqrt{\sigma^2+4\alpha D_S}}>0.
\end{equation*}
Moreover,
\begin{equation*}
\begin{aligned}
g_i(\chi_i^N-\chi_i^S)&=-\chi_i^S+\underbrace{\frac{\chi_i^N-\chi_i^S}{(\chi_i^N-\chi_i^S)^2+4\alpha D_S}}_{\leq 1}\chi_i^S<0,\\
g_i(\chi_i^N+\chi_i^S)&=\chi_i^S+\frac{\chi_i^N+\chi_i^S}{(\chi_i^N+\chi_i^S)^2+4\alpha D_S}\chi_i^S>0.
\end{aligned}
\end{equation*}
By a monotonicity argument, there exists a unique $\sigma_i \in I_i$ such that $g_i(\sigma_i)=0$. This corresponds to the individual speed of subpopulation $i$.\\
From the hypothesis that subpopulation 2 moves faster than subpopulation 1, we have $\sigma_1<\sigma_2$
From the dispersion relation \eqref{speed}, $g_1(\sigma)$ and $g_2(\sigma)$ are of opposite signs. Thus, from the monotonicity of $g_i$, we deduce that 
\begin{equation*}
\sigma \in (\sigma_1,\sigma_2).
\end{equation*}
Putting together the two conditions, we get
\begin{equation*}
\sigma \in (\sigma_1,\sigma_2)\cap (I_1\cap I_2).
\end{equation*}
If the two following conditions hold
\begin{equation*}
\begin{aligned}
\chi^N_2-\chi_2^S \not \in I_1\cap I_2\\
\sigma_2 \not \in I_1 \cap I_2 
\end{aligned}
\Rightarrow (\sigma_1,\sigma_2) \cap (I_1 \cap I_2)=(\sigma_1,\chi_1^N+\chi_1^S).
\end{equation*}
The dispersion relation \eqref{speed} is equivalent to
\begin{equation*}
G(\sigma)=\frac{\phi_{red}}{1-\phi_{red}},
\end{equation*}
where $G$ is defined by
\begin{equation*}
G(\sigma):=-\frac{\chi_2^S D_2}{\chi_1^S D_1}\frac{g_1(\sigma)}{g_2(\sigma)}\frac{h_2(\sigma)}{h_1(\sigma)}\frac{(\chi_1^S)^2-(\chi_1^N-\sigma)^2}{(\chi_2^S)^2-(\chi_2^N-\sigma)^2}.
\end{equation*}
Since $\sigma_2 \not \in (I_1 \cap I_2)$, we have $\chi^N_2+\chi_2^S \not \in (I_1\cap I_2)$ and the denominator of $G$ does not vanish
on $[\sigma_1,\chi_1^S+\chi^N_1]$. We conclude that the function $G$ is a bounded, positive and continuous function on $[\sigma_1,\chi_1^S+\chi^N_1]$.
Thus, $G$ attains its maximum over $[\sigma_1,\chi_1^S+\chi^N_1]$, denoted $\lambda^*$, and we have that
\begin{equation*}
\frac{\phi_{red}}{1-\phi_{red}} \leq \lambda^* \quad \Rightarrow \quad \phi_{red} \leq \phi_{red}^*=\frac{\lambda^*}{1+\lambda^*}.
\end{equation*}
This leads to the following conclusion
\begin{itemize}
\item $\displaystyle \phi_{red} \leq \phi_{red}^* \Rightarrow \quad $ Existence of a speed $\sigma$ satisfying the dispersion relation \eqref{speed}
\item $\displaystyle \phi_{red}>\phi_{red}^* \Rightarrow \quad$ Non-existence of a $\sigma$.
\end{itemize}

% Do NOT remove this, even if you are not including acknowledgments
\section*{Acknowledgments}
%Authors acknowledge financial support ANR-13-BS01-0004 Kibord from the French Ministry of Research and from the Labex IPGG. Ajouter INRIA
We thank V.Calvez, D.Lopez, J.Saragosti and P.Silberzan for numerous discussions and their help in setting up the experiments.
%\section*{References}
\bibliographystyle{unsrt}
\bibliography{plos_template}

\section*{Tables}
\begin{table}[!ht]
\caption{Parameter values}
\begin{tabular}{|c|c|c|}
\hline
Effective bacterial diffusion of subpopulation 1 $D_1$ &  $1.79 \times 10^{-6} \text{cm}^2.\text{s}^{-1}$ & Experimental measurement\\
\hline
Effective bacterial diffusion of subpopulation 2 $D_2$ &  $3.29 \times 10^{-6} \text{cm}^2.\text{s}^{-1}$ & Experimental measurement\\
\hline
Effective bacterial chemosensivity of subpopulation 1 $\chi_1^S$ &  $6.49\times 10^{-5} \text{cm}.\text{s}^{-1}$ & Experimental fit\\
\hline
Effective bacterial chemosensitivity of subpopulation 2 $\chi_2^S$ &  $2.88 \times 10^{-4} \text{cm}.\text{s}^{-1}$ & Experimental fit\\
\hline
Effective bacterial chemosensitivity of subpopulation 1 $\chi_1^N$ &  $2.57\times 10^{-4} \text{cm}.\text{s}^{-1}$ & Experimental fit\\
\hline
Effective bacterial chemosensitivity of subpopulation 2 $\chi_2^N$ &  $4.74 \times 10^{-4} \text{cm}.\text{s}^{-1}$ & Experimental fit\\
\hline
Chemical degradation $\alpha$ & $ 5 \times 10^{-2} \text{s}^{-1}$ & Experimental fit\\
\hline
Chemoattractant diffusion $D_S$ & $8 \times 10^{-6} \text{cm}^2.\text{s}^{-1}$ & \cite{Berg}\\
\hline
Nutrient diffusion $D_N$ & $8 \times 10^{-6} \text{cm}^2.\text{s}^{-1}$ & \cite{Berg} \\
\hline
\end{tabular}
\begin{flushleft}
\end{flushleft}
\label{parameters}
\end{table}

\end{document}